\documentclass{article}
\usepackage[utf8]{inputenc}
\usepackage[usenames]{color}
\usepackage{amsmath}
\usepackage{epsfig,latexsym}
\usepackage{epsf}
\usepackage{amsmath,amsfonts,amstext,amssymb,amsbsy,amsopn,amsthm,eucal}
\usepackage{latexsym}
\usepackage{epsfig,latexsym}
\usepackage[usenames]{color}
\usepackage{dirtytalk}

\usepackage{indentfirst}
\DeclareMathOperator{\tr}{tr}
\DeclareMathOperator{\sgn}{sgn}

\title{Generating Ouroboros Polynomials and Ouroboros Matrices}
\author{Nathan Thomas Provost$\footnote{Student of Applied Mathematics and Statistics at Brown University. University Email: nathan\_provost@brown.edu}$}
\date{}

\begin{document}

\maketitle

\begin{center}
\textbf{Abstract}
\end{center}

\small Ouroboros functions have shown some interesting properties when subjected to conventional operations. The aim of this paper is to continue our investigation and prove some additional properties of these functions. Using algebraic methods, we demonstrate that a collection of second-order polynomials can be generated for any multivariable Ouroboros function of the form we have mentioned in previous works (\cite{1} \cite{2}). We then generalize this observation to higher-order polynomials using the properties of Ouroboros spaces and the results of some of our previously proven theorems. After discussing the generation of these polynomials, we conclude by constructing a matrix from them and provide a few comments on its structure and aesthetic, culminating in the derivation of an intuitive formula for the degree of the trace of the square cases of these matrices and the discussion of some future research prospects.

\section*{Introduction}

\normalsize

It has been shown that the Ouroboros functions (as discussed in \cite{1} \cite{2} \cite{3}) possess some intriguing properties, especially after being subjected to conventional operations like differentiation (as in \cite{3}). In this case, we will continue to examine linear Ouroboros functions of the form:
\[f(x_1,...,x_n)=f(\textbf{x})=\sum_{i=1}^{n}c_ix_i \ni \sum_{i=1}^{n}c_i=1\]
Furthermore, we showed in \cite{2} that $f$ is an Ouroboros function for $\mathbb{R}^n$, or rather:
\[f\in\textbf{\textit{O}}(\mathbb{R}^n)=\left\{f:\mathbb{R}^n\rightarrow B \mid f(f(\textbf{x}),...,f(\textbf{x}))=f(\textbf{x}), \ \forall \textbf{x}\in \mathbb{R}^n, \ \forall B\subseteq \mathbb{R}\right\}\]
Following the logic presented in \cite{3}, we can impose several iterations of this defining feature, such that $f(f(f(...f(x)...)))=f(x)$ in one dimension for example. Using this property, we can derive polynomials in terms of the constants $\textbf{c}=(c_1,...,c_n)$, which possess some aesthetic and interesting properties. We will refer to these expressions as \textbf{Ouroboros Polynomials}.
\newpage

\section*{Quadratic Ouroboros Polynomials}

We begin by examining an instructive example in two dimensions. Consider the real-valued function $f=f(x,y)=c_1x+c_2y$ for which $c_1+c_2=1$. Following the theorem proven in \cite{2}, $f\in\textbf{\textit{O}}(\mathbb{R}^2)$, which means that $f(f,f)=f$. Intuitively, we see that $f(f,f)=c_1(c_1x+c_2y)+c_2(c_1x+c_2y)=(c_1+c_2)(c_1x+c_2y)$. Since $c_1+c_2=1$, we could omit the $(c_1+c_2)$ term, but for the purpose of exploration, we will hold off from doing so. Now, we see that:
\[ f(f(x,y),f(x,y))= f(x,y) \ \therefore \ (c_1+c_2)(c_1x+c_2y)  =  (c_1x+c_2y)\]
By equating coefficients, we see that $(c_1+c_2)c_2=c_2$ and $(c_1+c_2)c_1=c_1$, which makes sense since $c_1+c_2=1$. However, we can expand these polynomials to get the following quadratic equations: $c_1^2+c_2c_1-c_1=0$ and $c_2^2+c_1c_2-c_2=0$. These can be simplified to $c_1^2+(c_2-1)c_1=0$ and $c_2^2+(c_1-1)c_2=0$. Expanding upon this simple example, we can prove the following proposition.\\

\textbf{Proposition}: Suppose we have the following function for $n\in\mathbb{N} \ni n\geq2$:
\[f(x_1,...,x_n)=f(\textbf{x})=\sum_{i=1}^{n}c_ix_i \ni \sum_{i=1}^{n}c_i=1\]
where $\textbf{x}=(x_1,...,x_n)$, and suppose $k\in\{1,...,n\}\subset\mathbb{N}$. Furthermore, let $S_k=\{1,...,n\}/\{k\}$, $\forall k\in\{1,...,n\}$. Then:
\[c_1^2+\left(\left[\sum_{i\in S_1}c_i\right]-1\right)c_1=0\]
\[\vdots\]
\[c_n^2+\left(\left[\sum_{i\in S_n}c_i\right]-1\right)c_n=0\]

\noindent\textbf{Proof}: Assume all of the necessary assumptions and conditions of the proposition are met. As shown in the theorem from our previous paper \cite{2}, $f(\textbf{x})\in\textbf{\textit{O}}(\mathbb{R}^n)$, which means $f(f(\textbf{x}),...,f(\textbf{x}))=f(\textbf{x})$. Algebraically, it holds that:
\[f(f(\textbf{x}),...,f(\textbf{x}))=c_1\left(\sum_{i=1}^{n}c_ix_i\right)+...+c_n\left(\sum_{i=1}^{n}c_ix_i\right)=\]
\[\left(\sum_{i=1}^{n}c_i\right)\left(\sum_{i=1}^{n}c_ix_i\right)\]
\newpage
\noindent Therefore, we can further observe that:
\[ f(f(\textbf{x}),...,f(\textbf{x}))  =  f(\textbf{x}) \ \therefore \ \left(\sum_{i=1}^{n}c_i\right)\left(\sum_{i=1}^{n}c_ix_i\right)  =  \left(\sum_{i=1}^{n}c_ix_i\right)\] Focusing on $c_k$, we can equate coefficients and observe that:
\[\left(\sum_{i=1}^{n}c_i\right)c_k=c_k \therefore \left(\sum_{i=1}^{n}c_i\right)c_k-c_k=0\]
Furthermore, we can make an additional decomposition and determine that:
\[\sum_{i=1}^{n}c_i=c_k+\sum_{i\in S_k}c_i \ni \left(\sum_{i=1}^{n}c_i\right)c_k-c_k=\left(c_k+\sum_{i\in S_k}c_i\right)c_k-c_k=0 \therefore\]
\[c_k^2+c_k\left(\sum_{i\in S_k}c_i\right)-c_k=0 \therefore c_k^2+\left(\left[\sum_{i\in S_k}c_i\right]-1\right)c_k=0\]
Since $k$ was arbitrarily chosen from $\{1,...,n\}$, this proposition holds. $\qedsymbol$\\

\noindent\textbf{Remark}: We can solve this quadratic polynomial for $c_k$ using the quadratic formula, which yields:
\[c_k=\frac{1}{2}\left[\left(1-\left[\sum_{i\in S_k}c_i\right]\right)\pm\sqrt{\left(\left[\sum_{i\in S_k}c_i\right]-1\right)^2}\right]=\]
\[\frac{1}{2}\left[ \ \left(1-\left[\sum_{i\in S_k}c_i\right]\right)\pm\left|\left[\sum_{i\in S_k}c_i\right]-1\right| \ \right]\]
In this case, the sign of the expression within the absolute value bars is unimportant since the solutions are always given by:
\[c_k=\left(1-\left[\sum_{i\in S_k}c_i\right]\right) \ \rm{or} \ \textit{c}_\textit{k}=0 \ \rm{in \ either \ case.}\]
This is consistent with our assumption that the constants add up to 1, which further justifies the existence of our derived polynomial. As we have seen, one quadratic polynomial for each variable is generated, which makes sense, since this formula holds for any $k\in\{1,...,n\}$. Essentially, this establishes a set of $n$ relationships between the coefficients of a general, linear Ouroboros function. The next logical step in this investigation would be to consider a larger number of multiplicative iterations of the sum of these coefficients. As we will show in the pages that follow, this process yields higher-order polynomials that are satisfied by these coefficients.

\section*{Higher-Order Ouroboros Polynomials}

As we mentioned, the essential, defining property of an Ouroboros function is that $f(f(\textbf{x}),...,f(\textbf{x}))=f(\textbf{x})$ for any $\textbf{x}\in\mathbb{R}^n$. In this manner, the composition of an Ouroboros function with itself is an idempotent operation. Therefore, we can continuously compose the function with itself over an arbitrary number of iterations, generating an infinite number of higher-order Ouroboros polynomials. To illustrates this, let us revisit the function $f=f(x,y)=c_1x+c_2y$ where $c_1+c_2=1$. We can apply the properties of the Ouroboros functions to investigate the expression $f(f(f,f),f(f,f))=f$. Though it may seem complicated at first, we can see that $f(f(f,f),f(f,f))=c_1f(f,f)+c_2f(f,f)=(c_1+c_2)f(f,f)$ $=(c_1+c_2)(c_1f+c_2f)=(c_1+c_2)^2f=(c_1+c_2)^2(c_1x+c_2y)$. From this, we have:
\[ f(f(f,f),f(f,f)) = f \therefore   (c_1+c_2)^2(c_1x+c_2y) = (c_1x+c_2y) \]
By equating coefficients, we similarly find two polynomials: $(c_1+c_2)^2c_1=c_1$ and $(c_1+c_2)^2c_2=c_2$. After rearranging these expressions and expanding condensed terms, we have $c_1^3+2c_1^2c_2+c_1c_2^2-c_1=0$ and $c_2^3+2c_2^2c_1+c_2c_1^2-c_2=0$. These can be simplified to $c_1^3+(c_2^2+2c_1c_2-1)c_1=0$ and $c_2^3+(c_1^2+2c_1c_2-1)c_2=0$. These expressions are generated from the general relationships given by $c_1^3+((c_1+c_2)^2-c_1^2-1)c_1=0$ and $c_2^3+((c_1+c_2)^2-c_2^2-1)c_2=0$. These forms serve as special cases of the following  generalized theorem.\\

\textbf{Theorem}: Suppose we have a general, linear Ouroboros function:
\[f=f(x_1,...,x_n)=f(\textbf{x})=\sum_{i=1}^{n}c_ix_i\ni \sum_{i=1}^{n}c_i=1\]
Then for any $m\in\mathbb{N}$, for any $k\in\{1,...,n\}\subset\mathbb{N}$, and for any $n\in\mathbb{N}$:
\[c_k^m+\left(\left(\left[\sum_{i=1}^{n}c_i\right]^{m-1}-c_k^{m-1}\right)-1\right)c_k=0.\]
This expression is generated by equating the coefficients of
\[ \left[\sum_{i=1}^nc_i\right]^{m-1}\left(\sum_{i=1}^{n}c_ix_i\right)= \left(\sum_{i=1}^{n}c_ix_i\right),\]
which is valid since $f\in\textbf{\textit{O}}(\mathbb{R}^n)$ as proven in \cite{2}.\\

\noindent\textbf{Proof}: Assume that $m\in\mathbb{N}$, $k\in\{1,...,n\}\subset\mathbb{N}$, and $n\in\mathbb{N}$. Assume also that $\textbf{c}=(c_1,...,c_n)$ is a collection of real-valued coefficients that meets the necessary requirements for $f$. Let $\textbf{x}=(x_1,...,x_n)$. Denote $f(f(\textbf{x}),...,f(\textbf{x}))$ by $f_1(\textbf{F})$, denote $f(f_1(\textbf{F}),...,f_1(\textbf{F}))$ by $f_2(\textbf{F})$, and so on, such that the $p^{th}$ iteration of this complete self-composition is given by $f(f_{p-1}(\textbf{F}),...,f_{p-1}(\textbf{F}))=f_p(\textbf{F})$. Since we know 
$f\in\textbf{\textit{O}}(\mathbb{R}^n)$ from \cite{2}, we know $f_{m-1}(\textbf{F})=f(\textbf{x}), \forall m\in\mathbb{N}$ and:
\[f_{m-1}(\textbf{F})=\left[\sum_{i=1}^{n}c_i\right]^{m-1}\left(\sum_{i=1}^{n}c_ix_i\right)=\left(\sum_{i=1}^{n}c_ix_i\right)=f(\textbf{x})\]
As before, let $S_k=\{1,...,n\}/\{k\}$, $\forall k\in\{1,...,n\}$, so that we can write:
\[\sum_{i=1}^{n}c_ix_i=c_kx_k +\sum_{i\in S_k}c_ix_i\]
Simplifying this expression yields:
\[ f_{m-1}(\textbf{F})  =  f(\textbf{x})  \longrightarrow  \left[\sum_{i=1}^{n}c_i\right]^{m-1}\left(\sum_{i=1}^{n}c_ix_i\right)  =  \left(\sum_{i=1}^{n}c_ix_i\right) \]
\[\therefore \left[\sum_{i=1}^{n}c_i\right]^{m-1}  \left(c_kx_k +\sum_{i\in S_k}c_ix_i\right) = \left(c_kx_k +\sum_{i\in S_k}c_ix_i\right)\]
By equating coefficients and focusing on the $c_k$ terms, we see that:
\[\left[\sum_{i=1}^{n}c_i\right]^{m-1}c_k=c_k \therefore \left[\sum_{i=1}^{n}c_i\right]^{m-1}c_k-c_k=c_k^m+\left[\sum_{i=1}^{n}c_i\right]^{m-1}c_k-c_k^m-c_k=0\]
Therefore, through factorization, we have:
\[c_k^m+\left(\left(\left[\sum_{i=1}^{n}c_i\right]^{m-1}-c_k^{m-1}\right)-1\right)c_k=0. \ \qedsymbol\]

\noindent\textbf{Remark}: This result makes algebraic sense as well. We can see that:
\[\sum_{i=1}^{n}c_i=1 \longrightarrow \left[\sum_{i=1}^{n}c_i\right]^{m-1}=1, \ \rm{which \ in \ turn \ means \ that}:\]
\[c_k^m+\left(\left(\left[\sum_{i=1}^{n}c_i\right]^{m-1}-c_k^{m-1}\right)-1\right)c_k=c_k^m+(1-1-c_k^{m-1})c_k=c_k^m-c_k^m=0.\]
However, the goal of this derivation is not to redundantly make use of the fact that these constants sum to one. Instead, this theorem is a generalization of a relationship between the constants of all linear Ouroboros functions. Consequently, it describes a method for generating an arbitrary number of polynomials over an arbitrary number of iterations, the likes of which are always solved by the coefficients of the linear Ouroboros functions. Nonetheless, since this generalization is the result of iterative multiplication by the sum of all of the constants, it is algebraically sound \textbf{\textit{only}} because the sum of the constants is 1. We also note that in the trivial case where $n=1$, $c_1$ must be equal to 1, which would also justify the theorem, as the polynomial would reduce to $c_1^m-c_1=1^m-1=1-1=0$ as expected.

\section*{Ouroboros Matrices}

For any $k\in\{1,...,n\}\subset\mathbb{N}$ and for any $m\in\mathbb{N}$ with $j\in\{1,...,m\}\subset\mathbb{N}$, given a collection of real constants $\textbf{c}=(c_1,...,c_n)$ let us define:
\[p(k,j)=c_k^{j+1}+\left(\left(\left[\sum_{i=1}^{n}c_i\right]^{j}-c_k^{j}\right)-1\right)c_k.\]
where we assume that the constants meet the standard requirements for a linear Ouroboros function, such that:
\[\sum_{i=1}^nc_i=1 \ \ \rm{and} \ \ \textit{f}(\textit{x}_1,...,\textit{x}_\textit{n})=\sum_{\textit{i}=1}^\textit{n}\textit{c}_\textit{i}\textit{x}_\textit{i}\]
The \textbf{Ouroboros Matrix of Polynomials} (or more succinctly the \textbf{Ouroboros Matrix}) of dimension \textit{n} after \textit{m} iterations (denoted more explicitly by $\mathcal{M}[\textbf{\textit{O}}_{m}(\mathbb{R}^n)]$ and abbreviated by $\mathcal{M}$) is defined as:
\[\mathcal{M}[\textbf{\textit{O}}_{m}(\mathbb{R}^n)]=\begin{bmatrix}
p(1,1) & \dots & p(1,m) \\
\vdots & \ddots &\vdots \\
p(n,1) & \dots & p(n,m)
\end{bmatrix}\]
We note that this matrix does not contain the linear, or first-order polynomials, since these expressions reduce to $c_k-c_k$ for any $k$, which provides no insight into the relationships between coefficients. A case of particular interest (the \textbf{square} Ouroboros matrix) occurs when $m=n$, since under this condition we can consider the trace of $\mathcal{M}$. Moreover, for $\mathcal{M}[\textbf{\textit{O}}_n(\mathbb{R}^n)]=\mathcal{M}_{n^2}$:
\[\tr{(\mathcal{M}_{n^2})}=\prod_{i=1}^{n}p(i,i)=\left[c_1^{2}+\left(\left(\left[\sum_{i=1}^{n}c_i\right]-c_1\right)-1\right)c_1\right]\dots\left[c_n^{n+1}+\left(\left(\left[\sum_{i=1}^{n}c_i\right]^{n}-c_n^{n}\right)-1\right)c_n\right]\]
Naturally, the coefficients satisfy this equation for any $n$. It follows from this definition that since $\deg{(\tr{(\mathcal{M}_{n^2})})}=\deg{(p(1,1))}+...+\deg{(p(n,n))}$: 
\[\deg{(\tr{(\mathcal{M}_{n^2})})}=\sum_{i=2}^{n+1}i=n+\sum_{i=1}^{n}i=n+\frac{n(n+1)}{2}=\frac{n^2+3n}{2},\] 
resulting from the sum of the degrees of each $p(i,i)$, $\forall i\in\mathbb{N}$, which is a basic property of polynomials, and the universally accepted fact that the sum of the simplest arithmetic progression (stopping at some $n\in\mathbb{N}$) is equal to $\frac{1}{2}n(n+1)$, which historically has been proven through mathematical induction (a classical example of this proof can be found in sources like \cite{4}). Collectively, our definition of an Ouroboros matrix leaves us with this simple formula for determining the degree of the trace of a square Ouroboros matrix:
\[\deg{(\tr{(\mathcal{M}_{n^2})})}=\frac{n^2+3n}{2}, \ \forall n\in\mathbb{N}.\]

Let $\sigma$ represent a permutation, where $\sigma(i)$ is the post-permutation value of $i$, and let $P_n$ denote the set of all permutations for the set $\{1,...,n\}\subset\mathbb{N}$. Adapting the general definition of determinants from classical sources (like \cite{5} and \cite{6}), the determinant of $\mathcal{M}_{n^2}$ is:
\[|\mathcal{M}_{n^2}|=\sum_{\forall \sigma \in P_n}\left[\sgn{(\sigma)}\prod_{i=1}^{n}p(\sigma(i),i)\right]\]
where the sign function ($\sgn$) is defined as $\sgn{(\sigma)}=-1^{I_{\sigma}}$, where $I_{\sigma}$ is the number of inversions present in the permutation $\sigma$. Furthermore, the eigenvalues ($\lambda_1,...,\lambda_n$) of $\mathcal{M}_{n^2}$ are implicitly given as the solutions (in terms of $c_1,...,c_n$) to the polynomial generated by $|\mathcal{M}_{n^2}-\lambda \textbf{I}|=0$, where \textbf{I} is an $n\times n$ identity matrix (again, in accordance with traditional definitions like those from \cite{5} and \cite{6}).

\section*{Conclusion}

Evidently, we have shown that the coefficients of general, linear Ouroboros functions (with two or more variables) must satisfy an infinite number of polynomials. This results from the fact that we can multiply a linear Ouroboros function by the sum of its coefficients infinitely, since this sum has a value of one. Furthermore, we used these polynomials to define a matrix that will be of particular interest in future research endeavors. We showed that the degree of the trace of the square case of this matrix has a particularly aesthetic value, since it can be expressed as a concise quadratic expression in terms of the matrix's length and width ($n$), rather than a tedious expression that cannot be quickly evaluated for large values of $n$. For instance, the degree of the trace of a $100\times100$ Ouroboros matrix is 5150, which can be computed quickly due to the utility of the formula for $\deg{(\tr{(\mathcal{M}_{n^2})})}$. Additionally, we gave some brief insight into the determinants of these matrices by applying the traditional, permutation-based formula for the determinant of an $n\times n$ matrix to the entries of each Ouroboros matrix (given by $p(k,j)$). We also mentioned the characteristic polynomial through which the eigenvalues of these matrices can be obtained. While it is certainly possible to solve for the trivial eigenvalue (when $n=1$), it becomes increasingly difficult to solve for these eigenvalues as $n$ grows larger. This is due to the fact that the degrees of these characteristic polynomials grow much faster as $n$ increases. These concerns and problems will be the center of future research endeavors regarding Ouroboros matrices.

\end{document}